%% file: 2001-5.tex
\input gtmacros

\input pictex
\input rlepsf

\input gtoutput

\lognumber{136}
\volumenumber{5}\papernumber{5}\volumeyear{2001}
\pagenumbers{127}{142}
\accepted{6 March 2001}
\proposed{Cameron Gordon}
\seconded{Joan Birman, Wolfgang Metzler}
\received{12 September 2000}
\published{6 March 2001}


\def\gen#1{\langle #1 \rangle}
\def\ngen#1{\langle \langle #1 \rangle \rangle}
\let\la\langle
\let\ra\rangle
\let\wat\widehat

\font\spec=cmtex10 scaled 1095 
\def\d{\hbox{\spec \char'017\kern 0.05em}} 
\def\a{\longrightarrow}

\def\ex{\mathop{\rm ex}}
\def\inv{^{-1}}
\def\co{\mskip 1.5mu\colon\thinspace}
\def\H{{\cal H}}
\def\o{{\cal O}}
\def\Z{\Bbb Z}

\reflist

\key{Co} 
{\bf M\,M Cohen}, {\it A Course in Simple-Homotopy Theory},
Graduate Texts in Mathenmatics 10, Springer--Verlag, New York (1973)

\key{Coh}
{\bf M\,M Cohen}, {\it Whitehead torsion, group extensions, and
Zeeman's conjecture in high dimensions}, Topology, 16 (1977) 79--88

\key{CL}
{\bf Marshall M Cohen}, {\bf Martin Lustig}, {\it The conjugacy
problem for Dehn twist automorphisms of free groups}, Comment.
Math. Helv. 74 (1999) 179--200

\key{FR}
{\bf Roger Fenn}, {\bf Colin Rourke}, {\it Klyachko's methods and the
solution of equations over torsion-free groups}, L'Enseignment
Math\'ematique, 42 (1996) 49--74

\key{FRo}
{\bf Roger Fenn}, {\bf Colin Rourke}, {\it Characterisation of a class
of equations with solutions over torsion-free groups}, from: ``The
Epstein Birthday Schrift'', (I Rivin, C Rourke and C Series, editors),
Geometry and Topology Monographs, Volume 1 (1998) 163--17

\key{Me}
{\bf W Metzler}, {\it \"Uber den Homotopietyp zweidimensionaler
$CW$--Komplexe und Elementartransformationen bei Darstellungen von
Gruppen durch Erzeugende und definierende Relationen}, J. Reine
Angew. Math. 285 (1976) 7--23

\key{Ke}
{\bf M Kervaire}, {\it On higher dimensional knots}, from: ``Differential
and combinatorial top\-ology -- a symposium in honour of Marston Morse'',
Princeton Math. Series, 27 (1965)

\key{Kl}
{\bf A Klyachko}, {\it Funny property of sphere and equations over
groups}, Comm. in Alg. 21 (1993) 2555--2575

\key{Rot}
{\bf O\,S Rothaus}, {\it On the non-trivialty of some group extensions
given by generators and relations}, Ann. of Math. 106 (1977) 599--612

\key{RS}
{\bf C\,P Rourke}, {\bf B\,J Sanderson}, {\it Introduction to
piecewise-linear topology}, Sprin\-ger study edition, Springer--Verlag,
Berlin (1982)

\endreflist

\title{The surjectivity problem for one-generator,\\
one-relator extensions of torsion-free groups}
\shorttitle{Surjectivity problem for extensions of torsion-free groups}
\authors{Marshall M Cohen\\Colin Rourke}
\asciiaddress{Cornell University, Ithaca NY, 14853-4102, 
USA\\Warwick University,
Coventry, CV4 7AL, UK}
\address{Cornell University, Ithaca NY, 14853-4102, 
USA\\{\rm and}\\Warwick University,
Coventry, CV4 7AL, UK}
\email{marshall@math.cornell.edu\\cpr@maths.warwick.ac.uk}

\abstract
We prove that the natural map $G\to\wat G$, where $G$ is a
torsion-free group and $\wat G$ is obtained by adding a new generator
$t$ and a new relator $w$, is surjective only if $w$ is conjugate to
$gt$ or $gt\inv$ where $g\in G$.  This solves a special case of the
surjectivity problem for group extensions, raised by Cohen [\Coh].
\endabstract

\asciiabstract{%
We use Klyachko's methods to prove that the natural map G to G-hat,
where G is a torsion-free group and G-hat is obtained by adding a new
generator t and a new relator w, is surjective only if w is conjugate
to gt or gt^{-1} for some g in G.  This solves a special case of the
surjectivity problem for group extensions, raised by Cohen [Whitehead
torsion, group extensions, and Zeeman's conjecture in high dimensions,
Topology, 16 (1977) 79--88].}

\primaryclass{20E22, 20F05}\secondaryclass{57M20, 57Q10}
\asciikeywords{Surjectivity problem, torsion-free groups, Whitehead torsion,
Kervaire conjecture}

\maketitle

\section{Introduction}

In this paper we prove the following theorem.

\proclaim{Main Theorem}
Suppose that $G$ is a torsion-free group and that $\gen{t}$ is an
infinite cyclic group with generator t. Let $w$ be an
element of the free product  $G*\gen{t} $  and let $\ngen{w}$
be the normal subgroup of $G*\gen{t}$ generated by w.  View $G$
as a subgroup of $G*\gen{t}$ and let $i$ be the inclusion $G \a G*\gen{t}$.
Consider the natural homomorphism
$$q = \pi i\co \ G \ {\buildrel i \over \a} \ G*\gen{t}\ \
{\buildrel\pi\over \a}\ \ \wat{G} = {G*\gen{t}\over \ngen{w}}.$$ If
$q$ is onto then $w$ is conjugate to $gt$ or $gt\inv$ for some $g \in
G$.
\endproc

There are standard ways in which this algebraic situation may be
realized topologically. These lead to the following results.

\proclaim{Corollary 1} 
Suppose that $L$ is a connected CW complex with torsion-free
fundamental group and that the CW complex ${\wat L} = L \cup e^1 \cup
e^2$ is constructed by attaching a $1$--cell to $L$ and a 2--cell to
$L \cup e^1$.  If the inclusion map $j\co L \a \wat L$ induces a
surjection $j_*\co \pi_1L \a \pi_1{\wat L}$ then $j$ is a
simple-homotopy equivalence. \endproc

\prf 
This follows from elementary facts about the invariance of Whitehead
torsion under homotopy of attaching maps [\Co; Section 5] and the fact
that if ${\wat L} = L \cup e^1 \cup e^2$, where $e^1$ is a circle and
$e^2$ is a $2$-cell attached by a word $gt$, then $e^1$ is a free face
of $\wat L$ and $\wat L$ collapses to $L$ by an elementary
collapse. \endprf

\proclaim{Corollary 2} 
Suppose that $M$ is a connected $n$--manifold with torsion free
fundamental group and that $(W, M, M')$ is an h--cobordism with exactly
one handle of index one and one handle of index two and no other
handles (or dually with exactly one $n$--handle and one $ (n -
1)$--handle). Then $(W, M, M')$ is an s--cobordism.\endproc

\prf 
This is a consequence of the fact that a $k$-handle $D^k \times D^{n +
1 - k}$ collapses to its core union its attaching tube,
$D^k\times\{0\} \cup \partial D^k \times D^{n + 1 - k}$, see eg [\RS;
Chapter 6]. So the CW theory applies to the handlebody theory.\endprf

\sh{Background}

The surjectivity problem for group extensions and the question of
which Whitehead torsions can be realized were formulated by Cohen
[\Coh] and Metzler [\Me]; for more details on these problems and the
relevance of our results, see Section 5.

It will be useful to note from the outset that the conclusion
of the main theorem may be restated according to the following
lemma.

\proclaim{Lemma 1}If\ $G  \subset\, G*\gen{t}\ \
{\buildrel\pi\over \a}\ \ G*\gen{t}/\,\ngen{W}$\
where $W$ is a set of words in $G*\gen{t}$ then
$q = \pi\,\vert\,G$ is onto $\iff$ \ $gt$
lies in the kernel of $\pi$ for some $g \in G$.
\endproc
\prf $q$ is onto $\iff [\pi(t) \in \pi(G)] \iff [\pi(t) = \pi(g\inv)
\ {\rm for\ some}\ g \in G] \iff [\pi(gt) = 1] \iff
[gt \in {\rm kernel}(\pi)$\ for  some $g \in G]$. \endprf

Any element $w \in G*\gen{t}$ has a unique expression as a {\sl
reduced word}, $w = g_0t^{q_1}g_1t^{q_2}\ldots g_{n-1}t^{q_n}g_n$,
where $g_i\in G$ are non-trivial for $0<i<n$ and $q_i$ are non-zero
integers for each $i$.
The word $w$ is {\sl cyclically reduced} if further $g_n=1$ and if
$n>1$ then $g_0\neq1$.  Up to cyclic permutation there is a unique
cyclically reduced word in the conjugacy class of $w$, see eg [\CL,
Proposition 3.9].  
Since $\wat G$ depends only the conjugacy class of $w$, there is no
loss in assuming that $w$ is cyclically reduced and we shall do so
without comment from now on.  We call $\Sigma_{i=1}^nq_i$ the {\sl
exponent sum} of $t$ in $w$, denoted $\ex(w)$. The unreduced word
$t^{q_1}t^{q_2}\ldots t^{q_n}$ is called the {\sl $t$--shape} of $w$
and, thinking of $w=1$ as an equation over $G$, we call the elements
$g_i$ the {\sl coefficients} of $w$.

It is easy to see that if $q\co G \a \wat{G}$ is surjective then
$\ex(w) = \pm1$, since otherwise the abelianization of $\wat{G}/(q(G)
= 1)$ will be non-trivial.  So, replacing $w$ by $w^{-1}$ if
necessary, we may assume in our discussion that $\ex(w) = 1$.  Under
this hypothesis Klyachko [\Kl], in 1993, gave a brilliant argument to
prove the following theorem, which implies the Kervaire conjecture
[\Ke] in the case where $G$ is torsion-free.

\proclaim{Theorem}{\rm(Klyachko)}\qua 
If $G$ is a torsion-free group and $w \in G*\gen{t}$ with $\ex(w) = 1$
then the natural homomorphism $q\co G \a \wat{G} = {{G*\gen{t}}\over
\ngen{w}}$ is injective.
\endproc

An exposition (and extension) of Klyachko's theorem was given by Fenn
and Rourke [\FR] in 1996.  To prove our theorem we will use Klyachko's
result and his method, following closely the exposition in [\FR].  We
will quote some definitions and results from [\FR] and give those
proofs in detail for which the arguments differ and for which (proving
the contrapositive) the hypothesis is used that $w$ is not conjugate
to $gt$ for any $g \in G$.

\sh{Outline of the paper}

In Section 2 we consider a group $\Gamma$ in a slightly more general
situation than $G$ above. We assume (contrary to our Main Theorem)
that $w$ is not conjugate to $gt$ for any $g\in \Gamma$ but that some
$gt$ is in the kernel of $\Gamma*\gen{t} \a \wat \Gamma$. We show how
to construct a certain non-trivial CW subdivision of the $2$--sphere,
with edges labelled $t^{\pm1}$ and all but one corner labelled by an
element of $\Gamma$.

In Section 3 we prove our main theorem in a special case: We denote
$g^t = t\inv g t$. If $w$ has the form $w = b_0 a_0^t b_1 a_1^t \ldots
b_r a_r^t ct$, where the $a_i, b_i$ and $c$ are all elements of $G$
and $b_0 a_0^t b_1 a_1^t \ldots b_r a_r^t c \notin G$ then $q\co G \a
\wat{G}$ is not onto.

In Section 4, we complete the proof of the main theorem.  We use an
algebraic trick to parlay the result of Section 4 into a proof that,
in general, if $\ex(w) = 1$ and $w$ is not conjugate to $gt$ for any
$g \in G$ then $q\co G \a \wat{G}$ is not onto.

In Section 5 we briefly discuss the general surjectivity problem, in
which $n$ generators and $n$ relators are added to a group $G$.  We
give a bit of history and comment on the relevance of our result for
$n = 1$ to the general problem.

Finally, in Section 6 we extend our result to prove that if $w\in
G*\gen t$ is a word not of the form $gt$ whose $t$--shape is amenable
(see [\FR, \FRo]) then no word with $t$--shape $t^n$ can be in the
kernel of $\pi\co G\to\wat G = {G*\gen t\over\ngen w}$. 

\section{The cell subdivision lemma}

In this section we prove the cell subdivision lemma (below) which is
modelled on [\FR; Lemma 3.2].

The lemma uses the idea of a {\sl corner} of a 2--cell in a cell
subdivision $K$ of the 2--sphere.  This can be regarded as the
(oriented) angle formed by the two adjacent edges meeting at a vertex
(0--cell) in the boundary of the 2--cell.  If all the corners of a
2--cell are labelled by elements of a group, then a word can be read
around the 2--cell boundary by composing these elements either
unchanged or inverted according as the orientation of the corner
agrees or dissagrees with that of the 2--cell boundary.  Similarly if
all the corners at a vertex are labelled then a word can be read
around that vertex.  We shall always orient corners {\it clockwise},
thus if the above words are read {\it clockwise} for vertices and {\it
anticlockwise} for 2--cells, then no inversion is necessary.  See
figure 1 for an example: the word read around the boundary is $abc$;
after insertion of $t$ or $t\inv$ at the arrows (see part (e) of the
lemma below) it reads $tat\inv bt\inv c$.

\figure
\beginpicture
\setplotarea x from -2 to 2 , y from -0.5 to 4
\setcoordinatesystem units < .8500cm, .85000cm> point at 0 0
\linethickness=.5pt
\small
\putrule from -2 0 to 2 0
\setlinear
\plot -2 0  0 2  2 0 /
\circulararc 90 degrees from -.1 1.9 center at 0 2
\circulararc 45 degrees from -1.8 0 center at -2 0
\circulararc -45 degrees from 1.8 0 center at 2 0
\put {$\phantom a$} <0pt,0pt> [t] at 0 2  
\put {$a$} <0pt,-2pt> [t] at 0 1.8
\put {$b$} <1.5pt,0pt> [lb] at -1.7 0.1
\put {$c$} <-1.5pt,.5pt> [rb] at 1.7 0.1
\arrow <7pt> [.1,.5] from  -1 1 to -.9 1.1
\arrow <7pt> [.1,.5] from   1.1 .9 to 1 1
\arrow <7pt> [.1,.5] from .1 0 to -.1 0
\endpicture
\caption{Reading the boundary of a 2--cell: $tat\inv bt\inv c$}
\endfigure

Let $H$ be a subgroup of a group $\Gamma$ and let $g\in \Gamma$.  We
say that $g$ is {\sl free relative to} $H$ if the subgroup $\gen{g,H}$
of $\Gamma$ generated by $g$ and $H$ is naturally the free product
$\gen{g} * H$ of an infinite cyclic group $\gen g$ with $H$. (Note in
particular that $g$ has infinite order.)

If $g, h$ are elements of a group let $g^h$ denote $h^{-1}gh$.

In this section and the next, we shall consider the following
working hypotheses: 

\proclaim{Working hypotheses}

Suppose that $H$ and $H'$ are two isomorphic subgroups of a group
$\Gamma$ under the isomorphism $h\to h^{\phi}$, $h\in H$.  Suppose that
for each $i$, $a_i, b_i$ are elements of $\Gamma$ such that $a_i$ is
free relative to $H$ and $b_i$ is free relative to $H'$.  Let $c$ be
an arbitary element of $\Gamma$.

Let $w_0$ be the word $$b_0a_0^tb_1a_1^tb_2a_2^t\cdots b_ra_r^tct$$ in
$\Gamma*\gen t$, where $r\ge 0$, and and let $W\subset \Gamma*\gen t$
be the set of words $\{w_0,h^t(h^{\phi})\inv\mid h\in H\}$.  Let $\ngen W$
be the normal closure of $W$ in $\Gamma*\gen t$ and let $\wat\Gamma$
denote $(\Gamma*\gen{t})/\ngen W$.\endproc

\proclaim{Cell subdivision lemma}

Assume the working hypotheses, above.  Suppose that, for some $g\in
\Gamma$, $gt$ is in the kernel of the natural map $\Gamma*\gen
t\to\wat\Gamma$.

Then there is a cell subdivision $K$ of the 2--sphere such that

\items
\item{\rm(a)} the edges (1--cells) of $K$ are oriented,

\item{\rm(b)} the corners (all oriented clockwise) are labelled 
by coefficients of elements $w$ or $w\inv$ for $w\in W$, with the
exception of one particular corner at one particular vertex $v_0$
which is unlabelled,

\item{\rm(c)} the clockwise product of the corner labelling around any
vertex is $1\in\Gamma$ except for $v_0$ where it is undefined,

\item{\rm(d)} there is one special 2--cell $e^2_\infty$ which contains 
the unlabelled corner and has boundary a single edge and the single
vertex $v_0$,

\item{\rm(e)} with the exception of $e^2_\infty$, the corner labels of 
any 2--cell (in anticlockwise order) are the coefficients of $w$ or
$w\inv$ for some $w\in W$ (up to cyclic rotation) with the property
that, if on passing from one corner to an adjacent corner the element
$t$ or $t^{-1}$ is inserted according to whether the intervening edge
is oriented in the same or opposite direction (see figure 1), then the
whole of $w$ or $w\inv$ is recovered,

\item{\rm(f)} the cell decomposition is irreducible in the following
senses:\nl 
{\rm type (1)}\qua there do not exist two 2--cells with an edge in common
(necessarily read as $t$ in one and $t\inv$ in the other) such that,
starting with one vertex of this edge, the words read in these
2--cells are inverses of each other,\nl 
{\rm type (2)}\qua there does not exist a chain of 2--gons with common
vertices $a,b$ such that the product of the corner labels in
the chain at $a$ (or, equivalently, at $b$) is $1\in\Gamma$,

\item{\rm(g)} the cell subdivision is non-degenerate in that there
exist at least two vertices and at least three 2-cells; in particular
there is a cell $e^2_1 \neq e^2_\infty$ whose boundary contains
$\partial e^2_\infty$ as a {\it proper} subset (see figure 5).
\enditems\endproc

\prf  The proof uses transversality as in the proofs of [\FR; Lemmas 3.1 
and 3.2].

Choose a 2--complex $L$ with $\pi_1(L)=\Gamma$ and form the 2--complex
$\wat L$ with $\pi_1(\wat L)=\wat\Gamma$ by attaching a 1--cell
$\gamma$ to the base point $*$ of $L$ (corresponding to $t$) and a
2--cell $\sigma_w$ with attaching map determined by $w$ for each $w\in
W$.

Since $gt$ is trivial in $\pi_1(\wat L)$ there is a map of a 2--disc
$f\co D^2\to\wat L$ whose boundary maps to $L\cup\gamma$ and which
represents $gt\in\pi_1(L\cup\gamma)=\Gamma*\gen t$.  Make $f$
transverse to the centres of the 2--cells $\sigma_w$. It follows that
the inverse images of small neighbourhoods of these centres is a
collection of disjoint discs $D_1,\ldots,D_m$ in the interior of
$D^2$.  By a radial expansion of $f$ on these discs we may assume that
their image is the whole of one of the $\sigma_w$. It follows that the
punctured disc $P={\rm closure}\bigl(D^2- (D_1\cup\cdots\cup
D_m)\bigr)$ is mapped by $f$ to $L\cup\gamma$. Let $p$ be the centre
of $\gamma$. Make $f|P$ transverse to $p$. Then $f^{-1}p$ is a
1--manifold $Z$ properly embedded in $P$. By a radial expansion along
$\gamma$ we can assume that $Z$ has a neighbourhood $N$ which is a
normal $I$--bundle, where each fibre is mapped by $f$ to $\gamma$
and closure$(P - N)$ is mapped by $f$ to $L$.

We now simplify the subset $\H_f= D_1\cup\cdots\cup D_m\cup N$ of
$D^2$ as follows. Suppose $N$ contains an annulus component $\cal A$
in the interior of $P$. Let $D'$ denote the interior disc of $D^2$
which bounds the interior boundary component of the annulus. Then
$D''=D'\cup\cal A$ is a sub disc of $D^2$ whose boundary is mapped to
$*$ by $f$.  We can then redefine $f$ so that $f(D'')=*$ leaving $f$
unchanged outside $D''$.

At this point $\H_f$ can be regarded as a collection of 0--handles
(the $D_i$) and 1--handles (the components of $N$) attached to the
0--handles and to $\d D^2$ (in fact there is precisely one 1-handle
attached to $\d D^2$ by one end) see figure 2. 
\figure\relabelbox\small
\epsfxsize 2truein\epsfbox{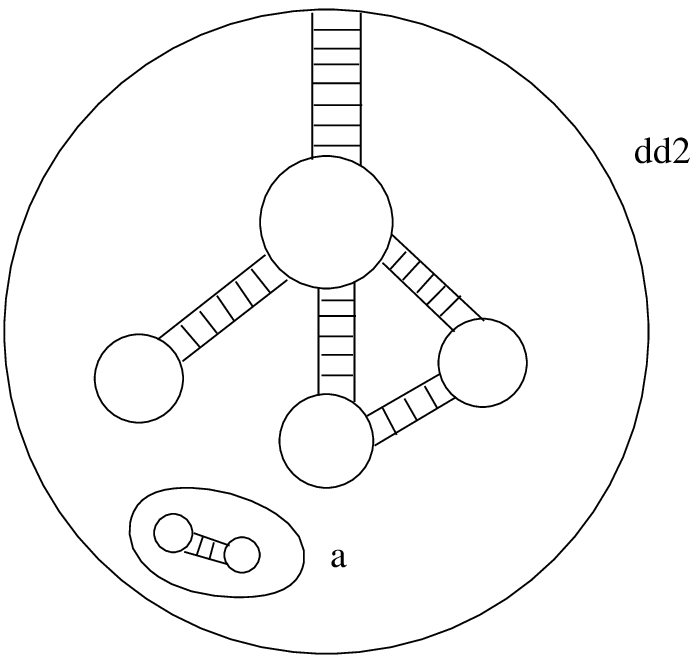}
\adjustrelabel <-4pt, 0pt> {a}{$\alpha$}
\adjustrelabel <-4pt, 0pt>  {dd2}{$\d D^2$}
\endrelabelbox
\caption {A view of $\H_f$}
\endfigure
We now prove that we
may assume that $\H_f$ is connected.  Suppose not.  Choose an
innermost component $C$.  Draw a simple loop $\alpha$ around $C$
separating it from the rest of $\H_f$.  Up to conjugacy $\alpha$
represents an element of $\pi_1(L)=\Gamma$ which is trivial in
$\pi_1(\wat L)=\wat\Gamma$.  But Klyachko proves that $\Gamma$ injects
in $\wat\Gamma$ (this is the precise content of [\FR; Theorem 4.1,
page 62]) and hence we may redefine $f$ so that the inside of $\alpha$
is mapped to $L$, which simplifies $\H_f$.

Note that the 0--handles can be labelled by elements $w$ or $w\inv$
for $w\in W$ according to the corresponding 2--cell of $\wat L$ and
orientation.   We say that
$\H_f$ is {\sl type (1) reducible} if there is a pair of 0--handles
labelled by $w$ and $w\inv$ (the same $w$) and joined by a 1-handle
which represents the same occurrence of $t$ (respectively $t\inv$) in
each word.  In this situation we can again simplify $\H_f$ without
changing $f|\d D^2$ by redefining $f$ near these 0--handles and
joining 1--handle (see figure 3).
\figure\relabelbox\small\let\ss\scriptstyle
\epsfxsize 4.5truein\epsfbox{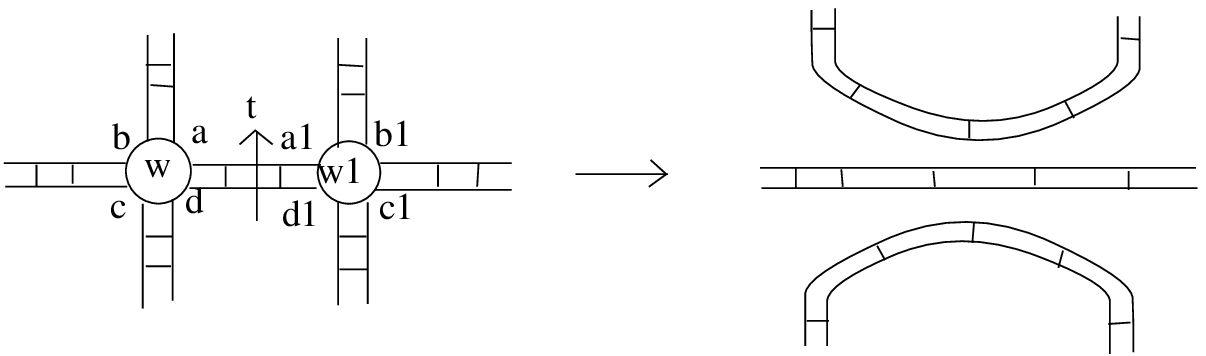}
\relabel {t}{$t$}
\relabel {a}{$\ss a$}
\relabel {b}{$\ss b$}
\relabel {c}{$\ss c$}
\relabel {d}{$\ss d$}
\relabel {w}{$\ss w$}
\relabel {a1}{$\ss a\inv$}
\relabel {b1}{$\ss b\inv$}
\relabel {c1}{$\ss c\inv$}
\relabel {d1}{$\ss d\inv$}
\relabel {w1}{$\ss w\inv$}
\endrelabelbox
\caption {Type (1) reduction of $\H_f$}
\endfigure
We say that $\H_f$ is {\sl type (2) reducible} if there is a chain of
0--handles (each having two 1--handles attached to it) labelled by
words $h_i^t(h_i^{\phi})\inv$, $i=1,2\ldots,q$ with $h_1h_2\ldots
h_q=1$ in $\Gamma$ (figure 4).
\figure\relabelbox\small\let\ss\scriptstyle
\epsfxsize 4.5truein\epsfbox{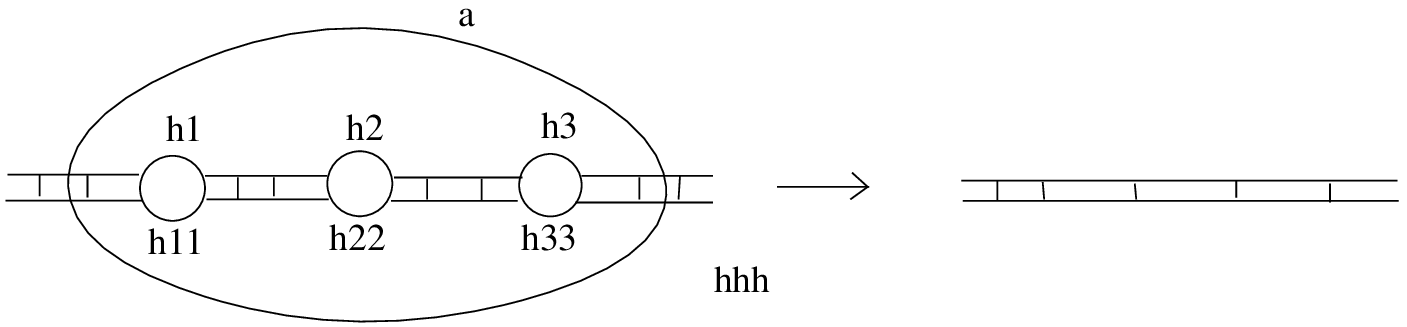}
\relabel {a}{$\alpha$}
\relabel {h1}{$\ss h_1$}
\relabel {h2}{$\ss h_2$}
\relabel {h3}{$\ss h_3$}
\relabel {h11}{$\ss (h_1^\phi)\inv$}
\relabel {h22}{$\ss (h_2^\phi)\inv$}
\relabel {h33}{$\ss (h_3^\phi)\inv$}
\relabel {hhh}{$h_1h_2h_3=1$}
\endrelabelbox
\caption {Type (2) reduction of $\H_f$}
\endfigure
If the chain forms a loop, the handlebody
is not connected and this chain and everything inside it may be
eliminated as indicated earlier.
Otherwise, the curve $\alpha$ indicated in figure 4 maps to
$tt\inv$ in $\Gamma*\gen t$ and there is another simplification given
by omitting this chain of 0--handles and redefining $f$ inside
$\alpha$ using the null-homotopy of $tt\inv$ in
$L\cup \gamma$. After these
simplifications there may now be more simplifications of the first two
types which can be performed.  Repeat all four until no more are
possible. Thus we can assume that $\H_f$ is connected and irreducible.

We now extend $\H_f$ to a handle decomposition $\H$ of $S^2\supset
D^2$ by letting the outside of $D^2$ be one 0--handle (denoted
$h^2_\infty$) and the regions of $D^2-\H_f$ be the 2--handles.

The required 2--complex $K$ is the dual complex to $\H$ obtained by
putting a vertex inside each 2--handle and joining by an edge across
each 1--handle.  The outside 2--cell is $e^2_\infty$ (containing
$h^2_\infty$) and has boundary containing a single vertex $v_0$.
Corners of 2--cells other than this corner are labelled by the
coefficient of the word $w$ or $w\inv$ labelling the 0--handle inside
the 2--cell opposite the corner.  See figure 5.

\figure\relabelbox\small\let\ss\scriptstyle
\epsfxsize 4truein\epsfbox{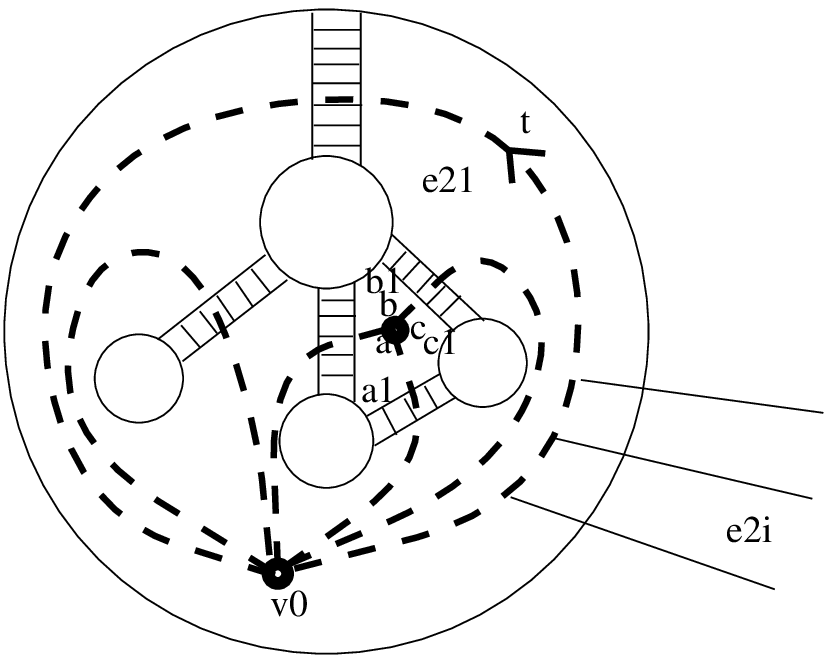}
\relabel {v0}{$v_0$}
\relabel {e21}{$e^2_1$}
\relabel {e2i}{$e^2_\infty$}
\relabel {t}{$t$}
\relabel {a}{$\ss a$}
\relabel {b}{$\ss b$}
\relabel {c}{$\ss c$}
\relabel {a1}{$\ss a$}
\relabel {b1}{$\ss b$}
\relabel {c1}{$\ss c$}
\endrelabelbox
\caption {The cell complex $K$ (shown dashed)}
\endfigure

The required properties of $K$ all follow from the construction:
1--cells are oriented by the orientations of the $I$--bundles
(1--handles) that they cross and properties (a) to (e) follow at once
(the word read around the boundary of a 2-cell is the label on the
contained 0--handle).  Property (f) follows from the irreducibility of
$\H_f$. 

Finally, property (g) uses the hypothesis that $r \geq 0$ (ie, that
$w_0$ is not conjugate to $gt$ for any $g \in \Gamma$). In order that
every $1$-handle of $\H$ have each end on some $\partial D^2_i$,
except that one of them has one end on $\partial D^2$, at least one of
the $0$-handles $D^2_i$ must be a $w_0$ or $w_0\inv$ handle. Because
$r\geq 0$, this handle must have at least three $1$-handles emanating
from it. Thus there have to be at least two $0$-handles inside $D^2$,
so that $K$ has at least three $2$-cells.  Since the handlebody closes
up, $D^2 - \H_f$ must have at least two components, resulting in at
least two vertices in $K$.
\endprf

\section{The key technical theorem}

In this section we prove the following result whose proof is modeled
on that of [\FR; Theorem 4.1].  We show that the hypotheses of the
cell subdivision lemma are self-contradictory.

\proclaim{Key Technical Theorem}

Assume the working hypotheses.  Then $gt$ is never in the kernel of
the natural map $\Gamma*\gen t\to\wat\Gamma$ for any $g\in \Gamma$.
\endproc

\rk{Remark}
Assuming this theorem, note that by Lemma 1 in the Introduction,
$\Gamma\to\wat\Gamma$ is not surjective. Therefore by taking $H$ and
$H'$ to be trivial, we can now deduce a special case of our main
theorem:
\smallskip
{\sl If the $t$--shape of $w = w_0$ is not $t$ (ie,
w is not conjugate to $gt$ for any $g \in \Gamma$) but is of the form
$t\inv tt\inv\ldots tt\inv tt$ then $q\co\Gamma\to\wat\Gamma$ 
is not surjective.}
\smallskip
In the next section we introduce an algebraic trick which will enable
us to deduce the general case, where $\ex(w)=1$ and $w$ is not
conjugate to $gt$ for any $g \in G$, from this special case.

\prf The proof relies heavily on the proof and terminology of
[\FR; Theorem 4.1, pages 62--64].
Assume that $gt$ is in the kernel of $\Gamma*\gen t\to\wat\Gamma$ where
$g\in\Gamma$.  By the cell subdivision lemma there is a cell subdivision
of $S^2$ with all 2-cells of the four types $I, I', II, II'$ illustrated
in Figure 6 with the exception of the special 2--cell $e^2_\infty$.  

\figure
\beginpicture
\setplotarea x from -2 to 5 , y from -4 to 4
\setcoordinatesystem units < 1.25000cm, 1.25000cm> point at 0 0
\linethickness=.5pt
\small
\putrule from -2 0 to 2 0
\put {$[2r+2,4r+1]$} [B] at 0 2.1  
\put {$0$} [t] at -2 -.1
\put {$1$} [t] at -1 -.1
\put {$2r$} [t] at 1 -.1
\put {$2r+1$} [t] at 2 -.1
\put {$c$} [t] at 0 1.8
\put {$b_0$} [lb] at -1.7 0.1
\put {$a_0$} [B] at -1 0.15
\put {$\ldots$} [B] at 0 0.15
\put {$a_r$} [rb] at 1.7 0.1
\put {$b_r$} [B] at 1 0.15 
\put {$I$} [] at 0 1
\setlinear
\plot -2 0  0 2  2 0 /
\circulararc 45 degrees from -1.8 0 center at -2 0
\circulararc -45 degrees from 1.8 0 center at 2 0
\circulararc 90 degrees from -.1 1.9 center at 0 2
\circulararc -180 degrees from -1.1 0 center at -1 0
\circulararc 180 degrees from 1.1 0 center at 1 0
\arrow <5pt> [.1,.5] from -1.5 0 to -1.6 0
\arrow <5pt> [.1,.5] from  -.5 0 to -.3 0
\arrow <5pt> [.1,.5] from  1.5 0 to 1.3 0
\arrow <5pt> [.1,.5] from  .4 0 to .6 0
\arrow <5pt> [.1,.5] from  -.9 1.1 to -1 1
\arrow <5pt> [.1,.5] from   1.1 .9 to 1 1
\setcoordinatesystem units < 1.3000cm, 1.3000cm> point at -5 -1
\circulararc -90 degrees from -1 0 center at 0 -1
\circulararc 90 degrees from -1 0 center at 0 1
\circulararc 90 degrees from .895 .095 center at 1 0
\circulararc 90 degrees from -.895 -.095 center at -1 0
\arrow <5pt> [.1,.5] from   0 .414 to -.141 0.407
\arrow <5pt> [.1,.5] from   0 -.414 to -.141 -0.407
\put {$II$} [] at 0 0
\put {$h^{-1}$} [] at .65 0 
\put {$h^{\phi}$} [] at -.7 0 
\put { } [] at 2 1
\setcoordinatesystem units < 1.25000cm, -1.25000cm> point at 0 -1
\putrule from -2 0 to 2 0
\put {$[1,2r]$} [t] at 0 2.1 
\put {$0$} [B] at -2 -.1
\put {$4r+1$} [B] at -1 -.1
\put {$2r+2$} [B] at 1 -.1
\put {$2r+1$} [B] at 2 -.1
\put {$c^{-1}$} [B] at 0 1.8
\put {$b_0^{-1}$} [lt] at -1.7 0.1
\put {$a_0^{-1}$} [t] at -1 0.15
\put {$\ldots$} [t] at 0 0.15
\put {$a_r^{-1}$} [rt] at 1.7 0.1
\put {$b_r^{-1}$} [t] at 1 0.15 
\put {$I'$} [] at 0 1
\setlinear
\plot -2 0  0 2  2 0 /
\circulararc -45 degrees from -1.8 0 center at -2 0
\circulararc 45 degrees from 1.8 0 center at 2 0
\circulararc -90 degrees from -.1 1.9 center at 0 2
\circulararc 180 degrees from -1.1 0 center at -1 0
\circulararc -180 degrees from 1.1 0 center at 1 0
\arrow <5pt> [.1,.5] from -1.5 0 to -1.6 0
\arrow <5pt> [.1,.5] from  -.5 0 to -.3 0
\arrow <5pt> [.1,.5] from  1.5 0 to 1.3 0
\arrow <5pt> [.1,.5] from  .4 0 to .6 0
\arrow <5pt> [.1,.5] from  -.9 1.1 to -1 1
\arrow <5pt> [.1,.5] from   1.1 .9 to 1 1
\setcoordinatesystem units < 1.3000cm, -1.3000cm> point at -5 -2
\circulararc 90 degrees from -1 0 center at 0 -1
\circulararc -90 degrees from -1 0 center at 0 1
\circulararc -90 degrees from .895 .095 center at 1 0
\circulararc -90 degrees from -.895 -.095 center at -1 0
\arrow <5pt> [.1,.5] from   0 .414 to -.141 0.407
\arrow <5pt> [.1,.5] from   0 -.414 to -.141 -0.407
\put {$II'$} <5pt,0pt> [] at 0 0
\put {$h^{}$} [] at .65 0 
\put {$({h^{\phi}})^{-1}$} <1pt,0pt> [] at -.5 0 
\put { } [] at 2 1
\endpicture
\caption{The 2-cells $I,I',II$ and $II'$}
\endfigure

As in [\FR] we give the two-sphere an orientation (``anticlockwise")
and give each $2$--cell of $K$ the induced orientation. A  traffic flow is
now defined, with a car running around the boundary of each $2$--cell in the
direction of the induced orientation  as follows:

At time 0 let a car on
the boundary of a country of type $I$ or $I'$ start at the corner labelled
$b_0$ or $b_0^{-1}$ and proceed in an anticlockwise manner with respect
to the orientation of the edge along which it is travelling, moving from
corner to corner in unit time except at the corner labelled $c$ or
$c^{-1}$ where it stops for $2r-1$ units. The times when the car is at
each corner are illustrated in figure 6. For countries of type $II$ or
$II'$ the car starts at the corner labelled $h^\phi$ or $({h^{\phi}})^{-1}$ 
and proceeds
in an anticlockwise manner moving from corner to corner in unit time.

For $e^2_\infty$ we need to consider also the 2--cell $e^2_1$ whose boundary
properly contains the boundary of $e^2_\infty$, see figure 5.  Let
$A$ be the car on the boundary of $e^2_\infty$ and $B$ the car on the
boundary of $e^2_1$.  Choose a point $\omega$ on $\d e^2_\infty$
different from $v_0$.  We shall engineer crashes between $A$ and $B$
to occur precisely at $\omega$.  Suppose that $B$ is approaching
$\omega$ then let $A$ approach $\omega$ from the opposite direction to
crash at $\omega$.  After the crash, let $A$ dawdle near $\omega$
until $B$ moves off $\d e^2_\infty$ (which it must, because $\d e^2_1$
{\it properly} contains $\d e^2_\infty$);  then let $A$ speed round to
just before $\omega$ where it again dawdles until $B$ again approaches
$\omega$ at which point the cycle repeats.

Recall from [\FR] that a {\it complete car crash} is said to occur when
two cars meet in the interior of an edge (necessarily going in opposite 
directions) or when a $N$ cars from $N$ neighbouring countries all meet
at a vertex of valency $N$.

Notice that, on $\d e^2_\infty$, the given flow has the property that
complete crashes occur at $\omega$ and nowhere else; in particular, no
complete crash occurs at $v_0$.  However there must be another
complete crash occurring at some other vertex of $K$. (This is for
exactly the same reasons as in the proof of [\FR; Theorem
4.1]. Properies 1 to 4 on pages 63--64 hold here also and the flow
satisfies the conditions of the Crash Theorem with Stops [\FR, Theorem
2.3, page 56].  So there must be another complete crash and, as in
[\FR], this must occur at a vertex.)

This leads to the identical contradiction as on [\FR; page 64]: The
flow has been chosen so that, at a vertex where all the cars come
together at the same time, the
labels around the corners are all $\{a,\, a\inv\}$ for some
coefficient $a=a_i$ of $w_0$ together with elements of $H$ or $\{b,\,
b\inv\}$ for some coefficient $b=b_i$ of $w_0$ together with elements
of $H'$.  For definiteness assume that we are in the former situation.
Then we can read an (unreduced) word of the form $a^\epsilon h_1 h_2
\dots h_{i_1} a^\epsilon h_1 h_2 \dots h_{i_2}a^\epsilon\dots$ which 
is 1 in $\Gamma$.  Now if this word contains a subword of the form
$a^\epsilon a^{-\epsilon}$ then $K$ is type (1) reducible and if it
contains a subword of the form $h_1 h_2 \dots h_i$ which is 1 in
$\Gamma$ then $K$ is type (2) reducible.  Since $K$ is irreducible
neither of these happen and the word either gives a non-trivial
relation in $\gen{a,H}$ contradicting the assumption that $a$ is free
relative to $H$ or reads $(a^\epsilon)^N = 1$ for $N\ge 1$ which also
contradicts the assumption that $a$ is free relative to $H$ (and in
particular has infinite order).
\endprf

\section{Proof of the main theorem}

In the light of the discussion in the Introduction and Lemma 1, we assume
that $\ex(w) = 1$ and, proving the contrapositive of the Main
Theorem, we assume
that $w$ is not conjugate to $gt$ for any $g \in G$. We must
prove that  
$$gt\hbox{\qua is not in the kernel of\qua }\pi\co G*\gen t\a \wat G
\ {\rm for \ any}\ g \in G. \eqno{(*)}$$

We now use Klyachko's algebraic trick described on pages 64--66 of
[\FR].  

Consider the homomorphism $\ex\co G*\gen t \to \Z$.  It is well known
that $K$, the kernel of $\ex$, is a free product of copies of $G$
generated by elements of the form $g^{t^\o}=t^{-\o}gt^\o$, $1 \neq g\in
G$.

Any element of $K$ has a canonical expression of the form
$k=g_1^{t^{\o_1}}\cdots g_r^{t^{\o_r}}$, where $\o_i\ne\o_{i+1}$ for
each $i$.   We shall call the $g_i^{t^{\o_i}}$ the {\sl canonical
elements} of $k$.  Let min$(k)$ be the minimum value of $\o_i$,
$i=1,\ldots,r$ and  max$(k)$ the maximum value. Fix a positive integer
$m$. Consider the following subgroups of $K$:
$$\eqalign{%
H&=\la k\in K\mid \,\hbox{\rm min}(k)\ge0, \hbox{\rm max}(k)\le m-2\ra  \cr
H'&=\la k\in K\mid \,\hbox{\rm min}(k)\ge1, \hbox{\rm max}(k)\le m-1\ra \cr
J&=\la k\in K\mid \,\hbox{\rm min}(k)\ge0, \hbox{\rm max}(k)\le m-1\ra  \cr}$$
and the following subsets:
$$\eqalign{%
X&=\{ k\in K\mid \,\hbox{\rm min}(k)=0, \hbox{\rm max}(k)\le m-1\}  \cr
Y&=\{ k\in K\mid \,\hbox{\rm min}(k)\ge0, \hbox{\rm max}(k)=m-1\}  \cr
Z&=\{ k\in K\mid \,\hbox{\rm min}(k)\ge1, \hbox{\rm max}(k)=m\}  . \cr}$$

\proclaim{Lemma 2}{\rm [\FR; Lemma 4.2, page 65]}\qua
Let $w\in G*\la t\ra $ satisfy $\ex(w)=1$. Then, after conjugation,
$w$ can be written as a product
$$b_0a_0^tb_1a_1^t\cdots b_ra_r^tct,$$
where $a_i\in Y, b_i\in X,i=0,\ldots, r$ and $c\in J$ for some $m>0$.

Furthermore, provided $w$ is not conjugate to $gt$ for some
$g\in G$, then $r\ge 0$ in the expression.\endproc 

\rk{Remark}The final sentence in the statement of lemma 2 is
not given in [\FR], but is immediate from the proof given there.

\proclaim{Lemma 3}{\rm [\FR; Lemma 4.3, page 66]}\qua
Suppose that $G$ is torsion-free, then any element $a$ of $Y$
is free relative to $H$.  Similary any element $b$ of $X$ is free
relative to $H'$.\qed\endproc

We can now complete the proof of the main theorem, which follows
closely the proof of [\FR; Theorem 4.4, pages 66--67].  By lemma 2 we
can assume that $w=b_0a_0^tb_1a_1^t\cdots b_ra_r^tct$, where $a_i\in
Y$, $b_i\in X$, $i=0,\ldots, r$ and $c\in J$ and $r\ge 0$.  We need to
think of each $a_i, b_i, c$ as functions of $t$ and for clarity we
shall introduce a new variable $s$.  To be precise let
$$w(s,t)\equiv b_0(t)a_0^s(t)\cdots b_r(t)a_r^s(t)c(t)s$$
where $s$ and $t$ are independent variables.

Write $\Gamma$ for $G*\la t\ra $ and let $H$, $H'$ be the subgroups
defined above.  There is an isomorphism $\phi:H\to H'$ given
by $h^\phi=h^t, h\in H$.

Lemma 3 gives the hypothesis of our key
technical theorem in Section 4, which implies that 
$$\hbox{$gs$\qua is never in the kernel of\qua $\Gamma*\gen s \to \wat 
\Gamma$.}\eqno{(**)}$$
The case $m=1$ is the special case (with $t$--shape $t\inv
tt\inv\ldots tt\inv tt$) covered by the proof in the last section, so
we may assume that $m>1$ and then $G\subset H\ne \emptyset$.

Each of the canonical elements of $a_i(t), b_i(t), c(t)$ is either in
$G$ or lies in $H^t$; moreover in $\wat{\Gamma}$ we have
$h^s=h^\phi=h^t$ for each $h\in H$.

Since $H$ is generated by elements of the form $t^{-i}gt^i$ for $i\le
m-2$ and since $h^s=h^t$ for each $h\in H$ it follows by induction on
$i$ that we can freely exchange $s$ and $t$ in products of elements of
the form $t^{-i}gt^i$ for $i\le m-1$. Thus we can exchange $s$ and $t$
in the coefficients of $w(s,t)$ and it follows that $w(s,s)=1$ in
$\wat{\Gamma}$.

Now consider the following commutative diagram 
$$\matrix{\Gamma=G*\gen t&\subset&\Gamma*\gen s&\a&\wat\Gamma\cr 
\cup&&\cup&&\uparrow\cr
G&\a&G*\gen s&\a&\widehat G\cr}$$ 
where $\widehat G={G*\gen s\over\ngen w}$.

By $(**)$, $gs\in \Gamma*\gen s$ does not map to $1\in \wat \Gamma$ 
for any $g \in G$. Therefore\ $gs \in G*\gen s$\ never maps to $1\in \wat G$.
This proves $(*)$ as required. \endprf

\section{The surjectivity problem and Whitehead torsion} 

If one adds $n$ generators $x_1, x_2, \ldots x_n$ and $n$ relators
$w_1, w_2, \ldots w_n$ to a group $G$ to form the group $\wat G$ then
one can ask whether the natural homomorphism $G \a \wat{G}$ is
injective.  If it is injective then one can ask whether it is
surjective.  Our main theorem answers this question completely for
torsion-free groups when $n = 1$.

The question of surjectivity, assuming injectivity, was raised by Cohen
[\Coh] in his study of Zeeman's conjecture.
Assuming the natural homomorphism is injective then one can associate
(after some normalization -- see [\Rot, pages 600-601]) to the set of
 words $w_1, \ldots w_n$ 
a Whitehead torsion element
$\tau \in {\rm Wh}(G)$ (the Whitehead group of $G$). Cohen conjectured
that if $\tau \neq 0$ then the injection cannot be onto.  In closely
related but independent work, in which he investigated 
inclusions of one 2-complex into another which are homotopy
equivalences, Metzler [\Me]  investigated the group theoretic
combinatorics and the set of Whitehead torsion elements which are associated to
such homotopy equivelences.  He named the set of Whitehead torsion elements
which can be realized by a relative 2--complex as ${\rm Wh}^\ast(G)$

Our main theorem appears to give evidence that  ${\rm Wh}^\ast(G) = 0$.
In fact, when $n = 1$ not only do we show that a necessary condition
for surjectivity is that the torsion of the $1 \times 1$ matrix is $0$,
but we show that, up to homotopy of the attaching map, the added one-
and two- cells can be collapsed away. However, one must be very cautious
concerning what this means for $n > 1$ in that

\items
\itemb not all Whitehead torsion
elements can be realized by $1\,\times \,1$ matrices, hence our result for
$n = 1$ in no way answers the question of whether ${\rm Wh}^\ast(G) = 0$, 

\itemb it is possible (an open conjecture) that ${\rm Wh}(G) = 0$ for all 
torsion-free groups.
\enditems

The only significant results on the surjectivity problem which we know
of for $n > 1$ are those of Rothaus [\Rot]. He develops an obstruction to
the surjectivity of the map $G \a \wat{G}$ in terms of representations
of $G$ into compact connected Lie groups. His theory had the following
application for dihedral groups, whose Whitehead groups are known to be
non-trivial.

\proclaim{Theorem}{\rm [\Rot; Theorem 11]}\qua If ${\rm p} \geq 5$ is an 
odd prime and $G = D_{2p}$ is the dihedral group of order 2$p$ and $n$
is any positive integer then there exist non-trivial Whitehead torsion
elements such that every injective homomorphism 
$G \a \wat{G} = {G*\gen{t_1, \ldots t_n}\over w_1 \ldots w_n}$
realizing this Whitehead torsion element is non-surjective.\endproc

Beyond Rothaus' work, the surjectivity problem for $n > 1$ is an
open and fascinating question.

\section{An extension and a question}

The main theorem (in the equivalent form given by lemma 1) can be
extended:

\proclaim{Extension of Main Theorem}

Let $G$ be a torsion free group and let $w\in G*\gen t$ be word which
is not of the form $gt$ and whose $t$--shape is amenable (see [\FR,
\FRo]) and consider the natural map $\pi\co G\to\wat G = {G*\gen
t\over\ngen w}$.  Let $x\in G*\gen t$ be any word with $t$--shape
$t^n$ for some $n>0$.

Then $x$ is not in the kernel of $\pi$.\endproc

The proof is very similar to the proof of the main theorem.  The cell
$e^2_\infty$ has $n$ edges all oriented the same way (``uphill'').
Notice that any other cell with an edge in common with $e^2_\infty$
has its car traverse that edge in the ``downhill" direction, since
adjacent cells induce opposite orientations on a common edge.  Choose
any point $\omega\in\d e^2_\infty$ not at a vertex.  The flow
constructed as in [\FR; page 68] for cells other than $e^2_\infty$ has
the property that there are times when all cars are going uphill and
hence are not on $\d e^2_\infty$.  This leaves time for car $A$ (on
$\d e^2_\infty$) to rush round from just after $\omega$ to just before
and hence there are no complete crashes on $\d e^2_\infty$ except at
$\omega$.  This leads to the identical contradiction as in the proof
of the main theorem.

The extension implies that all words of $t$--shape $t^n$ for some $n$ have
infinite order in $\wat G$.  This leads to the natural question:

\rk{Question}Suppose that $G$ is torsion-free and that $w$ is an
amenable word.  Is $\wat G$ torsion-free?
\ppar

If the answer is yes, then we can deduce that $G\to\wat G$ is never
surjective when $\wat G$ is obtained from $G$ by adding $n$ generators
and $n$ relators {\it one pair at a time}.

\references
\bye

%% file: gtmacros.tex
%
%
%
%
%
%
\magnification=\magstephalf      
%
%
\vsize=7.5truein                 
\hsize=5.2truein                 
\newskip\stdskip                 
\stdskip=6pt plus3pt minus3pt    
\medskipamount=\stdskip          
\parindent=0pt                   
\parskip=\stdskip                
\abovedisplayskip=\stdskip       
\belowdisplayskip=\stdskip       
\mathsurround=0.75pt             
\overfullrule=0pt                
%
%
\def\ppar{\par\goodbreak\vskip 8pt plus 4pt minus 4pt}     
%
%
\def\stdspace{\hskip 0.75em plus 0.15em\ignorespaces}
\let\qua\stdspace 
%
%
%
%
%
%
%
\def\hexnumber#1{\ifcase#1 0\or 1\or 2\or 3\or 4\or 5\or 6\or 7\or 8\or
 9\or A\or B\or C\or D\or E\or F\fi}
%
%
\font\thirtnmsa=msam10 scaled 1315    
\font\tenmsa=msam10          \font\ninemsa=msam9
\font\sevenmsa=msam7         \font\sixmsa=msam6
\font\fivemsa=msam5
%
%
\newfam\msafam                  \textfont\msafam=\tenmsa
\scriptfont\msafam=\sevenmsa    \scriptscriptfont\msafam=\fivemsa
\edef\hexa{\hexnumber\msafam}        
\def\msa{\fam\msafam\tenmsa}         
%
%
\font\thirtnmsb=msbm10 scaled 1315   
\font\tenmsb=msbm10      \font\ninemsb=msbm9
\font\sevenmsb=msbm7     \font\sixmsb=msbm6
\font\fivemsb=msbm5
%
\newfam\msbfam                   \textfont\msbfam=\tenmsb       
\scriptfont\msbfam=\sevenmsb     \scriptscriptfont\msbfam=\fivemsb
\edef\hexb{\hexnumber\msbfam}    
\def\msb{\fam\msbfam\tenmsb}     
%
%
\font\thirtneufm=eufm10 scaled 1315   
\font\teneufm=eufm10                 \font\nineeufm=eufm9
\font\seveneufm=eufm7                \font\sixeufm=eufm6
\font\fiveeufm=eufm5
%
\newfam\eufmfam                    \textfont\eufmfam=\teneufm
\scriptfont\eufmfam=\seveneufm     \scriptscriptfont\eufmfam=\fiveeufm
\edef\hexf{\hexnumber\eufmfam}      
\def\frak{\fam\eufmfam\teneufm}     
%
%
%
\font\thirtnrm=cmr10 scaled 1315    
\font\ninerm=cmr9                   \font\sixrm=cmr6   
%
\font\thirtni=cmmi10 scaled 1315    
\font\ninei=cmmi9                   \font\sixi=cmmi6  
%
\font\thirtnsy=cmsy10 scaled 1315   
\font\ninesy=cmsy9                  \font\sixsy=cmsy6  
%
\font\thirtnbf=cmbx10 scaled 1315   
\font\ninebf=cmbx9                  \font\sixbf=cmbx6  
%
%
\font\thirtnex=cmex10 scaled 1315   
\font\nineex=cmex9                  
%
%
\font\thirtnit=cmti10 scaled 1315  
\font\nineit=cmti9                  
%
\font\thirtnsl=cmsl10 scaled 1315  
\font\ninesl=cmsl9                  
%
\font\thirtntt=cmtt10 scaled 1315  
\font\ninett=cmtt9                  
%
%
%
%
\def\small{%
%
%
\textfont0=\ninerm \scriptfont0=\sixrm \scriptscriptfont0=\fiverm
\def\rm{\fam0\ninerm}
%
%
\textfont1=\ninei \scriptfont1=\sixi \scriptscriptfont1=\fivei
%
%
\textfont2=\ninesy \scriptfont2=\sixsy \scriptscriptfont2=\fivesy
%
%
\textfont3=\nineex \scriptfont3=\nineex \scriptscriptfont3=\nineex
%
%
\textfont\bffam=\ninebf \scriptfont\bffam=\sixbf
\scriptscriptfont\bffam=\fivebf \def\bf{\fam\bffam\ninebf}%
%
%
\textfont\itfam=\nineit \def\it{\fam\itfam\nineit}%
\textfont\slfam=\ninesl \def\sl{\fam\slfam\ninesl}%
\textfont\ttfam=\ninett \def\tt{\fam\ttfam\ninett}%
%
%
%
\textfont\msafam=\ninemsa \scriptfont\msafam=\sixmsa
\scriptscriptfont\msafam=\fivemsa \def\msa{\fam\msafam\ninemsa}%
%
%
\textfont\msbfam=\ninemsb \scriptfont\msbfam=\sixmsb
\scriptscriptfont\msbfam=\fivemsb \def\msb{\fam\msbfam\ninemsb}%
%
%
\textfont\eufmfam=\nineeufm  \scriptfont\eufmfam=\sixeufm
\scriptscriptfont\eufmfam=\fiveeufm \def\frak{\fam\eufmfam\nineeufm}%
%
%
%
\normalbaselineskip=11pt%
\setbox\strutbox=\hbox{\vrule height8pt depth3pt width0pt}%
%
%
\normalbaselines\rm
%
%
\stdskip=4pt plus2pt minus2pt    
\medskipamount=\stdskip          
\parskip=\stdskip                
\abovedisplayskip=\stdskip       
\belowdisplayskip=\stdskip       
\def\ppar{\par\goodbreak\vskip 6pt plus 3pt minus 3pt}%
%
%
\def\section##1{\global\advance\sectionnumber by 1
\vskip-\lastskip\penalty-800\vskip 20pt plus10pt minus5pt 
\egroup{\bf\number\sectionnumber\quad##1}\bgroup\small         
\vskip 6pt plus3pt minus3pt
\nobreak\resultnumber=1}
}    
%
\def\beginsmall{\bgroup\small}
\let\endsmall\egroup
%
%
%
%
\def\large{%
\textfont0=\thirtnrm \scriptfont0=\ninerm \scriptscriptfont0=\sevenrm
\def\rm{\fam0\thirtnrm}%
\textfont1=\thirtni \scriptfont1=\ninei \scriptscriptfont1=\seveni
\textfont2=\thirtnsy \scriptfont2=\ninesy \scriptscriptfont2=\sevensy
\textfont3=\thirtnex \scriptfont3=\thirtnex \scriptscriptfont3=\thirtnex
\textfont\bffam=\thirtnbf \scriptfont\bffam=\ninebf
\scriptscriptfont\bffam=\sevenbf \def\bf{\fam\bffam\thirtnbf}%
\textfont\itfam=\thirtnit \def\it{\fam\itfam\thirtnit}%
\textfont\slfam=\thirtnsl \def\sl{\fam\slfam\thirtnsl}%
\textfont\ttfam=\thirtntt \def\tt{\fam\ttfam\thirtntt}%
\textfont\msafam=\thirtnmsa \scriptfont\msafam=\ninemsa
\scriptscriptfont\msafam=\sevenmsa \def\msa{\fam\msafam\thirtnmsa}%
\textfont\msbfam=\thirtnmsb \scriptfont\msbfam=\ninemsb
\scriptscriptfont\msbfam=\sevenmsb \def\msb{\fam\msbfam\thirtnmsb}%
\textfont\eufmfam=\thirtneufm  \scriptfont\eufmfam=\nineeufm
\scriptscriptfont\eufmfam=\seveneufm \def\frak{\fam\eufmfam\teneufm}%
\normalbaselineskip=16pt%
\setbox\strutbox=\hbox{\vrule height11.5pt depth4.5pt width0pt}%
\normalbaselines\rm}%
\let\Large\large   
%
\def\Bbb#1{{\msb#1}}

%

%
\mathchardef\plussquare="0\hexa01
\mathchardef\nge="3\hexb0B
\mathchardef\maltesecross="0\hexa7A
\mathchardef\del="0\hexf01
%
%
%
%
\font\sc=cmcsc10
%
%
%
%
\def\sqr#1#2{{\vcenter{\vbox{\hrule  height.#2truept
	\hbox{\vrule width.#2truept height#1truept 
	\kern#1truept \vrule width.#2truept}
	\hrule height.#2truept}}}}
\def\sq{\sqr55}    
%
%
%
%
\newcount\sectionnumber            
\newcount\resultnumber             
\sectionnumber=0\resultnumber=1    
%
%
%
\def\section#1{\global\advance\sectionnumber by 1
\xdef\nextkey{\number\sectionnumber}
\vskip-\lastskip\penalty-800\vskip 20pt plus10pt minus5pt 
{\large\bf\number\sectionnumber\quad#1}         
\vskip 8pt plus4pt minus4pt
\nobreak\resultnumber=1}      
%
%
%
%
%
\def\sh#1{\vskip-\lastskip\ppar{\bf #1}\par\nobreak\medskip}         
%
%
%
%

%
\def\proc#1{\xdef\nextkey{\number\sectionnumber.\number\resultnumber}%
\vskip-\lastskip\ppar\bf%
\noindent#1\ \number\sectionnumber.\number\resultnumber
\stdspace\sl\global\advance\resultnumber by 1\ignorespaces}
\def\endproc{\rm\ppar} 
%
%
\def\prf{\vskip-\lastskip\ppar\noindent{\bf Proof}%
\stdspace\rm}                            
\def\qed{\hfill$\sq$\par\goodbreak\rm}   
\def\endprf{\unskip\stdspace\hbox{}
\hfill$\sq$\par\medskip}                 
%
%
%
%
%
%
%
%
\def\proclaim#1{\vskip-\lastskip\ppar\bf%
\noindent#1\stdspace\sl\ignorespaces} 

%
%
%
%
\def\rk#1{\vskip-\lastskip\ppar{\bf #1}\stdspace\ignorespaces}                

%
%
%
%
%
%
\def\label{\xdef\nextkey{\number\sectionnumber.\number\resultnumber}%
\number\sectionnumber.\number\resultnumber
\global\advance\resultnumber by 1}
%
%
%
%
%
%
%
%
%
%
%
%
%
%
%
%
\newcount\refnumber              
\refnumber=1                     
\long\def\reflist#1\endreflist{%
\long\def\thereflist{#1}{\def\refkey##1##2\par{\xdef##1{\number\refnumber}%
\global\advance\refnumber by 1}%
\def\key##1##2\par{\expandafter\xdef%
\csname##1\endcsname{\number\refnumber}%
\global\advance\refnumber by 1}#1\par}}
\long\def\references{%
\penalty-800\vskip-\lastskip\vskip 15pt plus10pt minus5pt 
{\large\bf References}\ppar 
{\leftskip=25pt\frenchspacing    
\small\parskip=3pt plus2pt       
\def\refkey##1##2\par{\noindent  
\llap{[##1]\stdspace}\ignorespaces##2\par}         
\def\key##1##2\par{\noindent  
\llap{[\ref{##1}]\stdspace}\ignorespaces##2\par}  
\def\,{\thinspace}\thereflist\par}}
%
%
%
\newcount\footnotenumber         
\footnotenumber=1                
\def\fnote#1{\xdef\nextkey{\number\footnotenumber}%
{\small\ifnum\footnotenumber>9\parindent=14pt%
\else\parindent=10pt\fi\footnote{$^{\number\footnotenumber}$}%
{\hglue-5pt#1}\global\advance\footnotenumber by 1}}
%
%
%
%
%
%
%
\newcount\figurenumber          
\figurenumber=1                 
\def\caption#1{\xdef\nextkey{\number\figurenumber}%
\cl{\small Figure \number\figurenumber: #1}%
\global\advance\figurenumber by 1}
\def\figurelabel{\xdef\nextkey{\number\figurenumber}%
\cl{\small Figure \number\figurenumber}%
\global\advance\figurenumber by 1}
\long\def\figure#1\endfigure{{\xdef\nextkey{\number\figurenumber}%
\let\captiontext\relax\def\caption##1{\xdef\captiontext{##1}}%
\midinsert\cl{\ignorespaces#1\unskip\unskip\unskip\unskip}\vglue6pt\cl{\small 
Figure \number\figurenumber\ifx\captiontext\relax\else: \captiontext
\fi}\endinsert\global\advance\figurenumber by 1}}
%
%
%
%
%
%
%
\def\nextkey{??}   
%
\def\key#1{\expandafter\xdef\csname #1\endcsname{\nextkey}}
\def\ref#1{\expandafter\ifx\csname #1\endcsname\relax
\immediate\write16{Reference {#1} undefined}??\else
\csname #1\endcsname\fi}
%
%
%
%
%
%
%
\newread\gtinfile
\newwrite\gtreffile
\def\useforwardrefs{
\openin\gtinfile\jobname.ref
\ifeof\gtinfile
\closein\gtinfile
\immediate\write16{No file \jobname.ref}
\else
\closein\gtinfile
\input \jobname.ref
\fi
\immediate\openout\gtreffile \jobname.ref
%
%
\def\key##1{{\def\\{\noexpand}%
\expandafter\xdef\csname ##1\endcsname{\nextkey}%
\immediate\write\gtreffile{\\\expandafter\\\def\\\csname ##1\\\endcsname%
{\nextkey}}}}
%
%
\long\def\reflist##1\endreflist{%
\long\def\thereflist{##1}{\def\refkey####1####2\par{\xdef####1{%
\number\refnumber}{\def\\{\noexpand}\immediate\write\gtreffile
{\\\def\\####1{\number\refnumber}}}\global\advance\refnumber by 1}%
\def\key####1####2\par{\expandafter\xdef%
\csname####1\endcsname{\number\refnumber}%
{\def\\{\noexpand}\immediate\write\gtreffile
{\\\expandafter\\\def\\\csname ####1\\\endcsname{\number\refnumber}}}
\global\advance\refnumber by 1}##1\par}}
\long\def\biblio##1\endbiblio{\reflist##1\endreflist\references}%
%
%
\def\numkey##1{{\def\\{\noexpand}%
\xdef##1{\number\sectionnumber.\number\resultnumber}
\immediate\write\gtreffile{\\\def\\##1%
{\number\sectionnumber.\number\resultnumber}}}}
\def\seckey##1{{\def\\{\noexpand}\xdef##1{\number\sectionnumber}
\immediate\write\gtreffile{\\\def\\##1{\number\sectionnumber}}}}
\def\figkey##1{\xdef##1{\number\figurenumber}%
{\def\\{\noexpand}\immediate\write\gtreffile%
{\\\def\\##1{\number\figurenumber}}}
\number\figurenumber\global\advance\figurenumber by 1}
}   
%
%
%
%
\def\figkey#1{\xdef#1{\number\figurenumber}%
\number\figurenumber\global\advance\figurenumber by 1}
\def\fig#1#2\endfig{%
\midinsert\cl{#2}\vglue6pt\cl{\small Figure #1}\endinsert}
\def\newfig{\number\figurenumber\global\advance\figurenumber by 1}
\def\numkey#1{\xdef#1{\number\sectionnumber.\number\resultnumber}}
\def\seckey#1{\xdef#1{\number\sectionnumber}}
%
%
%
%
%
%
%
%
%
\def\verb{\catcode`\"=\active}       
\def\brev{\catcode`\"=12}            
\brev                                
\verb                                
{\obeyspaces\gdef {\ }}              
{\catcode`\`=\active\gdef`{\relax\lq}}
\def"{%
\begingroup\baselineskip=12pt\def\par{\leavevmode\endgraf}%
\tt\obeylines\obeyspaces\parskip=0pt\parindent=0pt%
\catcode`\$=12\catcode`\&=12\catcode`\^=12\catcode`\#=12%
\catcode`\_=12\catcode`\~=12%
\catcode`\{=12\catcode`\}=12\catcode`\%=12\catcode`\\=12%
\catcode`\`=\active\let"\endgroup}
\brev      
%
%
%
%
%
%
\def\items{\par\leftskip = 25pt}           
\def\enditems{\par\leftskip = 0pt}         
\def\item#1{\par\leavevmode\llap{#1\stdspace}%
\ignorespaces}                             
\def\itemb{\item{$\bullet$}}               
%
%

%
%
\def\co{\colon\thinspace}    
\def\np{\vfil\eject}         
\def\nl{\hfil\break}         
\def\cl{\centerline}         
\def\gt{{\mathsurround=0pt\it $\cal G\mskip-2mu$eometry \&\ 
$\cal T\!\!$opology}}        
\def\agt{{\mathsurround=0pt\it$\cal A\mskip-.7mu$lgebraic \&\ 
$\cal G\mskip-2mu$eometric $\cal T\!\!$opology}}  
%
%
%

%
%
%
%
%
\def\title#1{\def\thetitle{#1}}
\def\shorttitle#1{\def\theshorttitle{#1}}
\def\author#1{\edef\previousauthors{\theauthors}
 \ifx\theauthors\relax\def\theauthors{#1}\else
 \def\theauthors{\previousauthors\par#1}\fi}

\let\authors\author        
\def\address#1{\edef\previousaddresses{\theaddress}
 \ifx\theaddress\relax\def\theaddress{#1}\else
 \def\theaddress{\previousaddresses\par\vskip 2pt\par#1}\fi}
\def\secondaddress#1{\edef\previousaddresses{\theaddress}
 \ifx\theaddress\relax\def\theaddress{#1}\else
 \def\theaddress{\previousaddresses\par{\rm and}\par#1}\fi}   

\def\email#1{\edef\previousemails{\theemail}
 \ifx\theemail\relax\def\theemail{#1}\else
 \def\theemail{\previousemails\hskip 0.75em\relax#1}\fi}
\def\secondemail#1{\edef\previousemails{\theemail}
 \ifx\theemail\relax\def\theemail{#1}\else
 \def\theemail{\previousemails\hskip 0.75em{\rm and}\hskip 0.75em
 \relax#1}\fi}
\def\url#1{\edef\previousurls{\theurl}
 \ifx\theurl\relax\def\theurl{#1}\else
 \def\theurl{\previousurls\hskip 0.75em\relax#1}\fi}
\def\secondurl#1{\edef\previousurls{\theurl}
 \ifx\theurl\relax\def\theurl{#1}\else
 \def\theurl{\previousurls\hskip 0.75em{\rm and}\hskip 0.75em
 \relax#1}\fi}
\long\def\abstract#1\endabstract{\long\def\theabstract{#1}}
\def\primaryclass#1{\def\theprimaryclass{#1}}
\def\secondaryclass#1{\def\thesecondaryclass{#1}}
\def\keywords#1{\def\thekeywords{#1}}
%
%
\let\\\par\let\thetitle\relax\let\theshorttitle\relax
\let\theauthors\relax\let\theshortauthors\relax
\let\theaddress\relax\let\theshortaddress\relax
\let\theemail\relax\let\theurl\relax
\let\theabstract\relax\let\theprimaryclass\relax
\let\thesecondaryclass\relax\let\thekeywords\relax
%
%
%
%
\long\def\maketitlepage{    

\vglue 0.2truein   

%
{\parskip=0pt\leftskip 0pt plus 1fil\def\\{\par\smallskip}{\large
\bf\thetitle}\par\medskip}   

\vglue 0.15truein 

%
{\parskip=0pt\leftskip 0pt plus 1fil\def\\{\par}{\sc\theauthors}
\par\medskip}%
 
\vglue 0.1truein 

%
{\small\parskip=0pt
{\leftskip 0pt plus 1fil\def\\{\par}{\sl\theaddress}\par}
\ifx\theemail\relax\else  
\vglue 5pt \def\\{\stdspace{\rm and}\stdspace} 
\cl{Email:\stdspace\tt\theemail}\fi
\ifx\theurl\relax\else    
\vglue 5pt \def\\{\stdspace{\rm and}\stdspace} 
\cl{URL:\stdspace\tt\theurl}\fi\par}

\vglue 7pt 

{\bf Abstract}

\vglue 5pt

\theabstract

\vglue 7pt 

{\bf AMS Classification numbers}\quad Primary:\quad \theprimaryclass\par

Secondary:\quad \thesecondaryclass

\vglue 5pt 

{\bf Keywords:}\quad \thekeywords

\np  

}    
%
%
\long\def\makeshorttitle{    


%
{\parskip=0pt\leftskip 0pt plus 1fil\def\\{\par\smallskip}{\large
\bf\thetitle}\par\medskip}   

\vglue 0.05truein 

%
{\parskip=0pt\leftskip 0pt plus 1fil\def\\{\par}{\sc\theauthors}
\par\medskip}%
 
\vglue 0.03truein 

%
{\small\parskip=0pt
{\leftskip 0pt plus 1fil\def\\{\par}{\sl\ifx\theshortaddress\relax
\theaddress\else\theshortaddress\fi}\par}
\ifx\theemail\relax\else  
\vglue 5pt \def\\{\stdspace{\rm and}\stdspace} 
\cl{Email:\stdspace\tt\theemail}\fi
\ifx\theurl\relax\else    
\vglue 5pt \def\\{\stdspace{\rm and}\stdspace} 
\cl{URL:\stdspace\tt\theurl}\fi\par}

\vglue 10pt 


{\small\leftskip 25pt\rightskip 25pt{\bf Abstract}\stdspace\theabstract

{\bf AMS Classification}\stdspace\theprimaryclass
\ifx\thesecondaryclass\relax\else; \thesecondaryclass\fi\par
{\bf Keywords}\stdspace \thekeywords\par}
\vglue 7pt
}    
\let\maketitle\makeshorttitle        
%
%

\def\volumenumber#1{\def\thevolumenumber{#1}}
\def\volumeyear#1{\def\thevolumeyear{#1}}
\def\pagenumbers#1#2{\def\startpage{#1}\def\finishpage{#2}}
\def\published#1{\def\publishdate{#1}}
\def\received#1{\def\receiveddate{#1}}

\let\reviseddate\relax
\volumenumber{X}
\volumeyear{20XX}
\pagenumbers{1}{XXX}
\published{XX Xxxember 20XX}

\long\def\makeagttitle{   
\agt\hfill      
\hbox to 60truept{\vbox to 0pt{\vglue -14truept{\bf [Logo here]}\vss}\hss}
\break
{\small Volume \thevolumenumber\ (\thevolumeyear)
\startpage--\finishpage\nl
Published: \publishdate}

\vglue .2truein

{\parskip=0pt\leftskip 0pt plus 1fil\def\\{\par\smallskip}{\large
\bf\thetitle}\par\medskip}   
\vglue 0.05truein 

%
{\parskip=0pt\leftskip 0pt plus 1fil\def\\{\par}{\sc\theauthors}
\par\medskip}%
 
\vglue 0.03truein 


{\small\leftskip 25truept\rightskip 25truept{\bf Abstract}\stdspace\theabstract

{\bf AMS Classification}\stdspace\theprimaryclass
\ifx\thesecondaryclass\relax\else; \thesecondaryclass\fi\par
{\bf Keywords}\stdspace \thekeywords\par}\vglue 7truept

}   


\def\Addresses{\bigskip
{\small \parskip 0pt \leftskip 0pt \rightskip 0pt plus 1fil \def\\{\par}
\sl\theaddress\par\medskip \rm Email:\stdspace\tt\theemail\par
\ifx\theurl\relax\else\smallskip \rm URL:\stdspace\tt\theurl\par\fi}}

\def\agtart{
\hoffset 14truemm
\voffset 31truemm
\font\phead=cmsl9 scaled 950
\font\pnum=cmbx10 scaled 913
\font\pfoot=cmsl9 scaled 950
\headline{\vbox to 0pt{\vskip -4.5mm\line{\small\phead\ifnum
\count0=\startpage ISSN numbers are printed here
\hfill {\pnum\folio}\else\ifodd\count0\def\\{ }%
\ifx\theshorttitle\relax\thetitle\else\theshorttitle\fi\hfill{\pnum\folio}
\else\def\\{ and }{\pnum\folio}\hfill\ifx\theshortauthors\relax\theauthors
\else\theshortauthors\fi\fi\fi}\vss}}
\footline{\vbox to 0pt{\vglue 0mm\line{\small\pfoot\ifnum\count0=\startpage
Copyright declaration is printed here\hfill\else
\agt, Volume \thevolumenumber\ (\thevolumeyear)\hfill\fi}\vss}}
\let\maketitle\makeagttitle\let\makeshorttitle\makeagttitle}

%% file: rlepsf.tex
\input epsf
%
%
\def\relabelbox{%
  \hbox\bgroup%
}%
\def\endrelabelbox{%
\egroup%
}%
\def\relabel #1#2 {%
  \special{ps:/a {} def}%
  \smash{\rlap{#2}}%
}%
\def\adjustrelabel <#1,#2> #3#4 {%
  \special{ps:/a {} def}%
  \smash{\rlap{\kern #1 \raise #2\hbox{#4}}}%
}%
\def\extralabel <#1,#2> #3 {\smash{\rlap{\kern #1 \raise #2\hbox{#3}}}}%

%% file: gtoutput.tex

\def\ifplaintex{\expandafter\ifx\csname documentclass\endcsname\relax}


\ifplaintex 
\hoffset 14truemm
\voffset 31truemm
\else
\headsep 23pt
\footskip 35pt
\hoffset -4truemm
\voffset 12.5truemm
\fi

\expandafter\ifx\csname beginpicture\endcsname\relax
\expandafter\ifx\csname documentclass\endcsname\relax
\input pictex \else\font\fiverm=cmr5
\input prepictex \input pictex \input postpictex \fi\fi

\def\gt{{\mathsurround=0pt\it $\cal G\mskip-2mu$eometry \&\ 
$\cal T\!\!$opology}}        

\def\gtp{{\mathsurround=0pt\it $\cal G\mskip-2mu$eometry \&\ 
$\cal T\!\!$opology $\cal P\!$ublications}}  


\def\lognumber#1{\def\thelognumber{#1}}
\def\volumenumber#1{\def\thevolumenumber{#1}}
\def\papernumber#1{\def\thepapernumber{#1}}
\def\volumeyear#1{\def\thevolumeyear{#1}}

\def\pagenumbers#1#2{\def\startpage{#1}\def\finishpage{#2}}
\def\published#1{\def\publishdate{#1}}
\def\proposed#1{\def\theproposer{#1}}
\def\seconded#1{\def\theseconders{#1}}
\def\received#1{\def\receiveddate{#1}}

\def\accepted#1{\def\accepteddate{#1}}

\def\asciiaddress#1{\def\theasciiaddress{#1}}

\long\def\asciiabstract#1{\long\def\theasciiabstract{#1}}
\def\asciikeywords#1{\def\theasciikeywords{#1}}

\def\shorttitle#1{\def\theshorttitle{#1}}


\let\\\par\let\thelognumber\relax
\let\thevolumenumber\relax\let\thepapernumber\relax
\let\thevolumeyear\relax\let\thesamplenumber\relax\let\startpage\relax
\let\finishpage\relax\let\publishdate\relax\let\receiveddate\relax
\let\reviseddate\relax\let\accepteddate\relax\let\theasciititle\relax
\let\theasciiauthors\relax\let\theasciiaddress\relax
\let\theasciiabstract\relax\let\theasciikeywords\relax
\let\theasciiemail\relax\let\theshortauthors\relax\let\theshorttitle\relax

\long\def\maketitlep{   

\count0=\startpage

\gt\hfill      
\beginpicture
\setcoordinatesystem units <0.33truein, 0.33truein> point at 2.2 0.9
\setplotsymbol ({$\cal G$})
\plotsymbolspacing=9truept
\circulararc 315 degrees from 0 1 center at 0 0
\setplotsymbol ({$\cal T$})
\circulararc 315 degrees from 1 -1 center at 1 0
\endpicture
%
\break
{\small\ifx\thesamplenumber\relax 
Volume \else Sample
\fi\thevolumenumber\ (\thevolumeyear)
\startpage--\finishpage\nl
Published: \publishdate}
\vglue 0.5truein plus 0.4fil minus 0.1truein

{\parskip=0pt\leftskip 0pt plus 1fil\def\\{\par\smallskip}{\ifplaintex\large
\else\Large\fi\bf\thetitle}\par\medskip}   

\vglue 0pt plus 0.1fil 

{\parskip=0pt\leftskip 0pt plus 1fil\def\\{\par}{\sc\theauthors}
\par\medskip}

\vglue 0pt plus 0.1fil 

{\small\parskip=0pt\let\newline\\
{\leftskip 0pt plus 1fil\def\\{\par}{\sl\theaddress}\par}
\expandafter\ifx\theemail\relax    
\relax\else\vglue 5pt plus 0.02fil minus 2pt\def\\{\stdspace{\rm 
and}\stdspace} 
\cl{Email:\stdspace\tt\theemail}\fi
\ifx\theurl\relax                  
\relax\else\vglue 5pt plus 0.02fil minus 2pt\def\\{\stdspace{\rm 
and}\stdspace}
\cl{URL:\stdspace\tt\theurl}\fi\par}

\vglue 7pt plus 0.3fil minus 3pt

{\bf Abstract}
\vglue 5pt plus 0.1fil minus 2pt

\theabstract

\vglue 7pt plus 0.3fil minus 3pt

{\bf AMS Classification numbers}\quad Primary:\quad \theprimaryclass

Secondary:\quad \thesecondaryclass

\vglue 5pt plus 0.3fil minus 2pt

{\bf Keywords}\quad \thekeywords

\vglue 10pt plus 0.5fil minus 5pt

{\small  Proposed: \theproposer\hfill Received: \receiveddate\nl
Seconded: \theseconders\hfill 
\ifx\reviseddate\relax                         
Accepted: \accepteddate                        
\else
Revised: \reviseddate                          
\fi}
\eject
}       

\let\maketitlepage\maketitlep
\let\maketitle\maketitlepage


\font\phead=cmsl9 scaled 950
\font\lhead=cmsl9 scaled 1050
\font\pnum=cmbx10 scaled 913
\font\lnum=cmbx10 
\font\pfoot=cmsl9 scaled 950
\font\lfoot=cmsl9 scaled 1050
\ifplaintex
\headline{\vbox to 0pt{\vskip -4.5mm\line{\small\phead\ifnum
\count0=\startpage ISSN 1364-0380 (on line)
1465-3060 (printed) \hfill {\pnum\folio}\else\ifodd\count0\def\\{ }%
\ifx\theshorttitle\relax\thetitle\else\theshorttitle\fi\hfill{\pnum\folio}
\else\def\\{ and }{\pnum\folio}\hfill\ifx\theshortauthors\relax\theauthors
\else\theshortauthors\fi\fi\fi}\vss}}
\footline{\vbox to 0pt{\vglue 0mm\line{\small\pfoot\ifnum\count0=\startpage
\copyright\ \gtp\hfill\else
\gt, Volume \thevolumenumber\ (\thevolumeyear)\hfill\fi}\vss
}}
\else
\makeatletter
\def\@oddhead{{\small\lhead\ifnum\count0=\startpage ISSN 1364-0380 (on line)
1465-3060 (printed) \hfill {\lnum\number\count0}\else\ifodd\count0
\def\\{ }\ifx\theshorttitle\relax \thetitle \else\theshorttitle\fi\hfill
{\lnum\number\count0}\else\def\\{ and }{\lnum\number\count0}
\hfill\ifx\theshortauthors\relax 
\theauthors\else\theshortauthors\fi\fi\fi}}\def\@evenhead{\@oddhead}
\def\@oddfoot{\small\lfoot\ifnum\count0=\startpage\copyright\ \gtp\hfill\else
\gt, Volume \thevolumenumber\ (\thevolumeyear)\hfill\fi}
\def\@evenfoot{\@oddfoot}
\makeatother
\fi


\newwrite\gtoutfile
\long\gdef\makeheadfile{  
{\def\\{, }\def\s{ }
\immediate\openout\gtoutfile head.xxx
\immediate\write\gtoutfile{To: math@arxiv.org}
\immediate\write\gtoutfile{Subject: put or rep NNNNN:pppp}
\immediate\write\gtoutfile{--text follows this line--}
\immediate\write\gtoutfile{Proxy-for: \ifx\theasciiauthors\relax
\theauthors\else\theasciiauthors\fi\s<\ifx\theasciiemail\relax\theemail\else\theasciiemail\fi>}
\immediate\write\gtoutfile{\noexpand\\}
\immediate\write\gtoutfile{Authors: \ifx\theasciiauthors\relax
\theauthors\else\theasciiauthors\fi}
\immediate\write\gtoutfile{Title: \ifx\theasciititle\relax
\thetitle\else\theasciititle\fi}
\immediate\write\gtoutfile{Subj-class: GT or SG or MG etc}
\immediate\write\gtoutfile{MSC-class: \theprimaryclass\ifx\thesecondaryclass\relax\else, \thesecondaryclass\fi}
\immediate\write\gtoutfile{Journal-ref: Geom. Topol. \thevolumenumber
(\thevolumeyear) \startpage-\finishpage}
\immediate\write\gtoutfile{Comments: Published by Geometry and Topology at}
\immediate\write\gtoutfile{\s\s http://www.maths.warwick.ac.uk/gt/GTVol\thevolumenumber/paper\thepapernumber.abs.html}
\immediate\write\gtoutfile{\noexpand\\}
\immediate\write\gtoutfile{}
\ifx\theasciiabstract\relax
\immediate\write\gtoutfile{\theabstract}\else
\immediate\write\gtoutfile{\theasciiabstract}\fi
\immediate\write\gtoutfile{}
\immediate\write\gtoutfile{\noexpand\\}
\immediate\write\gtoutfile{}
\immediate\closeout\gtoutfile}}  

\def\maketitlepage{\maketitlep\makeheadfile}
\let\maketitle\maketitlepage


\def\ifplaintex{\expandafter\ifx\csname documentclass\endcsname\relax}


\ifplaintex 
\hoffset 14truemm
\voffset 31truemm
\else
\headsep 23pt
\footskip 35pt
\hoffset -4truemm
\voffset 12.5truemm
\fi

\expandafter\ifx\csname beginpicture\endcsname\relax
\expandafter\ifx\csname documentclass\endcsname\relax
\input pictex \else\font\fiverm=cmr5
\input prepictex \input pictex \input postpictex \fi\fi

\def\gt{{\mathsurround=0pt\it $\cal G\mskip-2mu$eometry \&\ 
$\cal T\!\!$opology}}        

\def\gtp{{\mathsurround=0pt\it $\cal G\mskip-2mu$eometry \&\ 
$\cal T\!\!$opology $\cal P\!$ublications}}  


\def\lognumber#1{\def\thelognumber{#1}}
\def\volumenumber#1{\def\thevolumenumber{#1}}
\def\papernumber#1{\def\thepapernumber{#1}}
\def\volumeyear#1{\def\thevolumeyear{#1}}

\def\pagenumbers#1#2{\def\startpage{#1}\def\finishpage{#2}}
\def\published#1{\def\publishdate{#1}}
\def\proposed#1{\def\theproposer{#1}}
\def\seconded#1{\def\theseconders{#1}}
\def\received#1{\def\receiveddate{#1}}

\def\accepted#1{\def\accepteddate{#1}}

\def\asciiaddress#1{\def\theasciiaddress{#1}}

\long\def\asciiabstract#1{\long\def\theasciiabstract{#1}}
\def\asciikeywords#1{\def\theasciikeywords{#1}}

\def\shorttitle#1{\def\theshorttitle{#1}}


\let\\\par\let\thelognumber\relax
\let\thevolumenumber\relax\let\thepapernumber\relax
\let\thevolumeyear\relax\let\thesamplenumber\relax\let\startpage\relax
\let\finishpage\relax\let\publishdate\relax\let\receiveddate\relax
\let\reviseddate\relax\let\accepteddate\relax\let\theasciititle\relax
\let\theasciiauthors\relax\let\theasciiaddress\relax
\let\theasciiabstract\relax\let\theasciikeywords\relax
\let\theasciiemail\relax\let\theshortauthors\relax\let\theshorttitle\relax

\long\def\maketitlep{   

\count0=\startpage

\gt\hfill      
\beginpicture
\setcoordinatesystem units <0.33truein, 0.33truein> point at 2.2 0.9
\setplotsymbol ({$\cal G$})
\plotsymbolspacing=9truept
\circulararc 315 degrees from 0 1 center at 0 0
\setplotsymbol ({$\cal T$})
\circulararc 315 degrees from 1 -1 center at 1 0
\endpicture
%
\break
{\small\ifx\thesamplenumber\relax 
Volume \else Sample
\fi\thevolumenumber\ (\thevolumeyear)
\startpage--\finishpage\nl
Published: \publishdate}
\vglue 0.5truein plus 0.4fil minus 0.1truein

{\parskip=0pt\leftskip 0pt plus 1fil\def\\{\par\smallskip}{\ifplaintex\large
\else\Large\fi\bf\thetitle}\par\medskip}   

\vglue 0pt plus 0.1fil 

{\parskip=0pt\leftskip 0pt plus 1fil\def\\{\par}{\sc\theauthors}
\par\medskip}

\vglue 0pt plus 0.1fil 

{\small\parskip=0pt\let\newline\\
{\leftskip 0pt plus 1fil\def\\{\par}{\sl\theaddress}\par}
\expandafter\ifx\theemail\relax    
\relax\else\vglue 5pt plus 0.02fil minus 2pt\def\\{\stdspace{\rm 
and}\stdspace} 
\cl{Email:\stdspace\tt\theemail}\fi
\ifx\theurl\relax                  
\relax\else\vglue 5pt plus 0.02fil minus 2pt\def\\{\stdspace{\rm 
and}\stdspace}
\cl{URL:\stdspace\tt\theurl}\fi\par}

\vglue 7pt plus 0.3fil minus 3pt

{\bf Abstract}
\vglue 5pt plus 0.1fil minus 2pt

\theabstract

\vglue 7pt plus 0.3fil minus 3pt

{\bf AMS Classification numbers}\quad Primary:\quad \theprimaryclass

Secondary:\quad \thesecondaryclass

\vglue 5pt plus 0.3fil minus 2pt

{\bf Keywords}\quad \thekeywords

\vglue 10pt plus 0.5fil minus 5pt

{\small  Proposed: \theproposer\hfill Received: \receiveddate\nl
Seconded: \theseconders\hfill 
\ifx\reviseddate\relax                         
Accepted: \accepteddate                        
\else
Revised: \reviseddate                          
\fi}
\eject
}       

\let\maketitlepage\maketitlep
\let\maketitle\maketitlepage


\font\phead=cmsl9 scaled 950
\font\lhead=cmsl9 scaled 1050
\font\pnum=cmbx10 scaled 913
\font\lnum=cmbx10 
\font\pfoot=cmsl9 scaled 950
\font\lfoot=cmsl9 scaled 1050
\ifplaintex
\headline{\vbox to 0pt{\vskip -4.5mm\line{\small\phead\ifnum
\count0=\startpage ISSN 1364-0380 (on line)
1465-3060 (printed) \hfill {\pnum\folio}\else\ifodd\count0\def\\{ }%
\ifx\theshorttitle\relax\thetitle\else\theshorttitle\fi\hfill{\pnum\folio}
\else\def\\{ and }{\pnum\folio}\hfill\ifx\theshortauthors\relax\theauthors
\else\theshortauthors\fi\fi\fi}\vss}}
\footline{\vbox to 0pt{\vglue 0mm\line{\small\pfoot\ifnum\count0=\startpage
\copyright\ \gtp\hfill\else
\gt, Volume \thevolumenumber\ (\thevolumeyear)\hfill\fi}\vss
}}
\else
\makeatletter
\def\@oddhead{{\small\lhead\ifnum\count0=\startpage ISSN 1364-0380 (on line)
1465-3060 (printed) \hfill {\lnum\number\count0}\else\ifodd\count0
\def\\{ }\ifx\theshorttitle\relax \thetitle \else\theshorttitle\fi\hfill
{\lnum\number\count0}\else\def\\{ and }{\lnum\number\count0}
\hfill\ifx\theshortauthors\relax 
\theauthors\else\theshortauthors\fi\fi\fi}}\def\@evenhead{\@oddhead}
\def\@oddfoot{\small\lfoot\ifnum\count0=\startpage\copyright\ \gtp\hfill\else
\gt, Volume \thevolumenumber\ (\thevolumeyear)\hfill\fi}
\def\@evenfoot{\@oddfoot}
\makeatother
\fi


\newwrite\gtoutfile
\long\gdef\makeheadfile{  
{\def\\{, }\def\s{ }
\immediate\openout\gtoutfile head.xxx
\immediate\write\gtoutfile{To: math@arxiv.org}
\immediate\write\gtoutfile{Subject: put or rep NNNNN:pppp}
\immediate\write\gtoutfile{--text follows this line--}
\immediate\write\gtoutfile{Proxy-for: \ifx\theasciiauthors\relax
\theauthors\else\theasciiauthors\fi\s<\ifx\theasciiemail\relax\theemail\else\theasciiemail\fi>}
\immediate\write\gtoutfile{\noexpand\\}
\immediate\write\gtoutfile{Authors: \ifx\theasciiauthors\relax
\theauthors\else\theasciiauthors\fi}
\immediate\write\gtoutfile{Title: \ifx\theasciititle\relax
\thetitle\else\theasciititle\fi}
\immediate\write\gtoutfile{Subj-class: GT or SG or MG etc}
\immediate\write\gtoutfile{MSC-class: \theprimaryclass\ifx\thesecondaryclass\relax\else, \thesecondaryclass\fi}
\immediate\write\gtoutfile{Journal-ref: Geom. Topol. \thevolumenumber
(\thevolumeyear) \startpage-\finishpage}
\immediate\write\gtoutfile{Comments: Published by Geometry and Topology at}
\immediate\write\gtoutfile{\s\s http://www.maths.warwick.ac.uk/gt/GTVol\thevolumenumber/paper\thepapernumber.abs.html}
\immediate\write\gtoutfile{\noexpand\\}
\immediate\write\gtoutfile{}
\ifx\theasciiabstract\relax
\immediate\write\gtoutfile{\theabstract}\else
\immediate\write\gtoutfile{\theasciiabstract}\fi
\immediate\write\gtoutfile{}
\immediate\write\gtoutfile{\noexpand\\}
\immediate\write\gtoutfile{}
\immediate\closeout\gtoutfile}}  

\def\maketitlepage{\maketitlep\makeheadfile}
\let\maketitle\maketitlepage